\documentclass[12pt]{article}
\usepackage{amsmath, amssymb, amsfonts}
\usepackage[utf8]{inputenc}
\usepackage{geometry}
\usepackage{hyperref}
\usepackage{url}
\title{Constructive Limits of Cantor's Diagonal Method: Countability, Enumerability, and the Impossibility of Exhausting the Continuum}
\author{Stanislav Semenov \\
\href{mailto:stas.semenov@gmail.com}{stas.semenov@gmail.com} \\
\href{https://orcid.org/0000-0002-5891-8119}{ORCID: 0000-0002-5891-8119}}
\date{March 23, 2025}

\begin{document}

\maketitle

\begin{abstract}
  Cantor's diagonal method is traditionally used to prove the uncountability of the set of all infinite binary sequences. This paper analyzes the expressive limits of this method. It is shown that under any constructive application—including generalizations with computable permutations and infinite hierarchies of diagonal extensions—the resulting set remains countable. Thus, the method demonstrates the incompleteness of countable coverage but is unable to generate an uncountable set. This highlights its limitations as a constructive tool and reveals the boundary between constructive enumerability and the completeness of the continuum.
\end{abstract}

\subsection*{Mathematics Subject Classification}
03D80 (Computability and recursion), 03E10 (Ordinal and cardinal numbers), 03B70 (Logic in computer science)

\subsection*{ACM Classification}
F.4.1 Mathematical Logic, F.1.1 Models of Computation

\section*{Introduction}

Cantor's diagonal method is widely known as one of the most elegant proofs of the existence of uncountable sets—in particular, the set of all infinite binary sequences \(\mathbb{B}^\mathbb{N}\)~\cite{cantor1891}. Its classical application demonstrates that no countable enumeration can exhaust the entire space \(\mathbb{B}^\mathbb{N}\), since it is always possible to construct an element not included in the given enumeration.

However, a natural question arises: what exactly does this method create—does it truly "generate" uncountability, or does it merely emphasize its unattainability through constructive processes? Can the diagonal method not only demonstrate incompleteness but also construct an uncountable set?

In this article, we explore the limits of the expressive power of Cantor's diagonal method. We formalize the set of all diagonally constructed elements, including those generalized using computable permutations, and trace how even under infinite hierarchical application, the method remains strictly within the framework of countable constructions. We show that the diagonal method alone cannot generate an uncountable set and, consequently, cannot be used as a criterion for the uncountability of arbitrary sets. In this sense, it serves as a tool for demonstration but not for construction.

\section{Construction of the Set \texorpdfstring{\(X\)}{X} and the Classical Diagonal Method}

Let
\[
X = \{x_1, x_2, x_3, \dots\} \subset \mathbb{B}^\mathbb{N}
\]
be a countable set, where each element is an infinite binary sequence:
\[
x_i = (a_{i1}, a_{i2}, a_{i3}, \dots), \quad a_{ij} \in \{0, 1\}
\]

Cantor's diagonal method constructs a new sequence
\[
y = (b_1, b_2, b_3, \dots), \quad b_i = 1 - a_{ii}
\]

Thus, \(y\) differs from each \(x_i\) at least in the \(i\)-th position, and therefore:
\[
y \notin X
\]

\section{Diagonal Counterexamples \texorpdfstring{\(Y\)}{Y}}

Define \(Y \subset \mathbb{B}^\mathbb{N}\) as the set of all sequences \(y_k\), each constructed using the diagonal principle but with different permutations\footnote{Only \emph{computable} permutations \(\pi: \mathbb{N} \to \mathbb{N}\) are considered, since for the constructive construction of \(y_k\) by the rule \(b_i^{(k)} = 1 - a_{i \pi(i)}\), it is necessary that the value \(\pi(i)\) can be computed for any \(i\). The set of all such computable permutations is countable (this follows from the fact that each computable permutation can be specified by an algorithm, and the set of all algorithms is countable~\cite{sipser, soare}), which ensures the countability of \(Y\). If arbitrary (including non-computable) bijections \(\mathbb{N} \to \mathbb{N}\) are allowed, the set \(Y\) may become uncountable, but then it loses constructive definiteness, and the diagonal method ceases to be applicable in an explicit form.} of index pairs \((i,j)\) in the matrix \(a_{ij}\). Formally, for each permutation \(\pi_k: \mathbb{N} \to \mathbb{N}\), define the sequence \(y_k\) as follows:
\[
y_k = (b_1^{(k)}, b_2^{(k)}, b_3^{(k)}, \dots), \quad b_i^{(k)} = 1 - a_{i \pi_k(i)}
\]

Examples:
\begin{itemize}
  \item \(y_1\): the main diagonal \(\pi_1(i) = i\), i.e., \(b_i^{(1)} = 1 - a_{ii}\)
  \item \(y_2\): the permutation \(\pi_2(1) = 2, \pi_2(2) = 1, \pi_2(i) = i\) for \(i \geq 3\), i.e., \(b_1^{(2)} = 1 - a_{12}, b_2^{(2)} = 1 - a_{21}, b_i^{(2)} = 1 - a_{ii}\) for \(i \geq 3\)
  \item \(y_3\): the permutation \(\pi_3(1) = 3, \pi_3(2) = 2, \pi_3(3) = 1, \pi_3(i) = i\) for \(i \geq 4\), i.e., \(b_1^{(3)} = 1 - a_{13}, b_2^{(3)} = 1 - a_{22}, b_3^{(3)} = 1 - a_{31}, b_i^{(3)} = 1 - a_{ii}\) for \(i \geq 4\)
\end{itemize}

Thus,
\[
Y = \{y_1, y_2, y_3, \dots\}
\]
is a countable set whose elements are constructively defined based on the elements of \(X\) and a predetermined diagonal permutation.

\textbf{Proof of the Countability of \(Y\)}: Since each sequence \(y_k\) is uniquely determined by the permutation \(\pi_k\), and the set of all permutations \(\pi_k\) is countable (as a subset of \(\mathbb{N} \to \mathbb{N}\)), \(Y\) is also countable.

\section{Application of the Diagonal Method to the Set \texorpdfstring{\(Y\)}{Y}}

Consider the set \(Y = \{y_1, y_2, y_3, \dots\}\), where each sequence \(y_i\) is constructed using the diagonal principle with some permutation \(\pi_i: \mathbb{N} \to \mathbb{N}\). That is:

\[
y_i = (b_1^{(i)}, b_2^{(i)}, b_3^{(i)}, \dots), \quad b_j^{(i)} = 1 - a_{j, \pi_i(j)}
\]

Each \(y_i\) differs from the row \(x_j\) at the index \(j = \pi_i^{-1}(i)\), if such a permutation is invertible, or at a specially constructed diagonal set of indices. At the same time, each \(y_i\) is uniquely and constructively determined by the corresponding permutation.

Apply the diagonal method to the set \(Y\), using as the "diagonal" the elements \(b_i^{(i)}\) from the \(i\)-th row \(y_i\) and position \(i\):

\[
z = (c_1, c_2, c_3, \dots), \quad c_i = 1 - b_i^{(i)} = a_{i, \pi_i(i)}
\]

Thus:

\begin{itemize}
  \item Each \(c_i\) is a value taken from the row \(x_i\) (an element of the set \(X\)) at position \(\pi_i(i)\)
  \item The resulting sequence \(z\) is assembled from individual bits, each of which belongs to some \(x_i \in X\), but the sequence as a whole is not required to coincide with any single \(x_j\)
\end{itemize}

For \(z \in X\) to hold, there must exist some \(k\) such that:
\[
\forall i \in \mathbb{N},\quad a_{i, \pi_i(i)} = a_{k,i}
\]
However, since each \(c_i\) can be taken from its own row \(x_i\), and the structure of the set \(X\) is arbitrary, such a coincidence is not guaranteed.

Therefore, in general:
\[
z \notin X
\]

The diagonal method, applied to the set \(Y\), yields a new sequence \(z\) that may not belong to the original set \(X\). Thus, even if all \(y_i\) were constructed from elements of \(X\), the result of their "reverse diagonal inversion" does not necessarily return us to \(X\).

\textbf{Remark.} The resulting sequence \(z\) also does not belong to \(Y\), since each of its elements is constructed in the opposite direction relative to the construction of all \(y_i\). That is:
\[
z \notin Y, \quad \text{and in general } z \notin X
\]
which confirms the ability of the diagonal method to generate new elements outside any given countable set, even if that set already contains all diagonally inverted forms of other elements.

\section{Iterative Application of the Diagonal Method and the Boundary of Constructivity}

Consider the union of the original set \(X\) and the set of all diagonally inverted sequences \(Y\), constructed using various computable permutations:

\[
X_1 := X \cup Y
\]

The set \(X_1\) remains countable, since the union of two countable sets is also countable. Apply the diagonal method to \(X_1\):

\[
w_1 := \text{diag}(X_1)
\]

Since \(w_1\) is constructed to differ from at least one element of each row in \(X_1\), it does not belong to \(X_1\):
\[
w_1 \notin X_1
\]

Add this new element:
\[
X_2 := X_1 \cup \{w_1\}
\]

Repeat the process: apply the diagonal method to \(X_2\), obtain \(w_2 := \text{diag}(X_2)\), which is not in \(X_2\), and construct:
\[
X_3 := X_2 \cup \{w_2\}
\]

Continuing by induction, define an ascending sequence of sets:
\[
\begin{aligned}
X_1 &:= X \cup Y \\
X_{n+1} &:= X_n \cup \{\text{diag}(X_n)\}
\end{aligned}
\]

Each \(X_n\) is countable and contains all previous sets: \(X_n \subsetneq X_{n+1}\), since the diagonal method guarantees an element outside the current set. All new elements \(w_n\) are constructively obtained and uniquely determined by \(X_n\)~\cite{odifreddi}.

Define the limit set as the union of all such iterative extensions:

\[
X_\infty := \bigcup_{n=1}^\infty X_n
\]

\textbf{Countability of \(X_\infty\)}: Since each \(X_n\) is countable, and the union of a countable number of countable sets is also countable, we have:
\[
|X_\infty| = \aleph_0
\]

\section{Limitations of the Method and the Impossibility of Exhausting \texorpdfstring{\(\mathbb{B}^\mathbb{N}\)}{B(N)}}

Despite the infinite iteration of diagonal extensions, the set \(X_\infty\) remains countable. However:

\[
|\mathbb{B}^\mathbb{N}| = 2^{\aleph_0} \gg \aleph_0 \Rightarrow X_\infty \subsetneq \mathbb{B}^\mathbb{N}
\]

Thus, no constructive procedure based on the diagonal method and its repeated applications can cover the entire set of infinite binary sequences.

\section{Enumerability and Non-Computability}

Each set \(X_n\), as well as the limit set \(X_\infty\), is constructively enumerable—its elements are obtained through an algorithmic process (diagonal construction and union). Therefore:

\[
X_\infty \subset \mathcal{R}
\]

where \(\mathcal{R}\) denotes the set of all computable sequences.

Since the set \(\mathcal{R}\) is countable, and the set \(\mathbb{B}^\mathbb{N}\) is of continuum cardinality, the inevitable conclusion is that the vast majority of sequences in \(\mathbb{B}^\mathbb{N}\) are non-computable and cannot be generated by any algorithm.

\textbf{Discussion of Non-Computable Elements.} Non-computable sequences are those for which no algorithm exists that can enumerate their bits. Their existence is guaranteed by the fact that the set of all programs (or Turing machines) is countable~\cite{sipser, davis}, while the set of all infinite binary sequences has the cardinality of the continuum~\cite{soare, pour_el}. Thus, most elements of \(\mathbb{B}^\mathbb{N}\) are fundamentally unreachable by means of constructive construction, including the diagonal method.

\section{Cantor's Method and the Limits of Its Expressiveness}

Cantor's diagonal method, when applied to a countable set of binary sequences, generates a new sequence guaranteed not to belong to the original set. This provides a powerful tool for proving the uncountability of the set \(\mathbb{B}^\mathbb{N}\) as a whole.

However, as shown above, the method does not transcend the constructive paradigm: it takes a countable set as input and returns a specific sequence, the addition of which preserves countability. Even an infinite hierarchy of extensions based on the sequential application of the diagonal method,
\[
X_1 \subset X_2 \subset X_3 \subset \dots \subset X_\infty = \bigcup_{n=1}^\infty X_n
\]
remains within the bounds of countable cardinality: \(|X_\infty| = \aleph_0\).

This leads to an important limitation:

\begin{quote}
\textbf{The diagonal method cannot generate an uncountable set.}
\end{quote}

This means that the diagonal method cannot:

\begin{itemize}
  \item construct an uncountable set through iterations;
  \item enumerate all elements of \(\mathbb{B}^\mathbb{N}\);
  \item be used as a test for the uncountability of an arbitrary set \(X \subset \mathbb{B}^\mathbb{N}\).
\end{itemize}

Even if the set \(X\) is not explicitly given but merely assumed as a candidate for uncountability, the application of the diagonal method does not provide a criterion for its verification. The method can show that a particular enumeration is \emph{incomplete}, but it cannot determine whether the set \(X\) itself is countable or uncountable.

\textbf{Conclusion.} Cantor's diagonal method possesses expressive power only within the bounds of countable constructions and cannot encompass the entirety of \(\mathbb{B}^\mathbb{N}\). Its limitations highlight the distinction between constructive and non-constructive aspects of infinity and underscore the need for other, more general tools for analyzing uncountable sets.

\section*{Conclusion}

Cantor's diagonal method is a powerful way to demonstrate that any countable enumeration of elements of the set \(\mathbb{B}^\mathbb{N}\) is inevitably incomplete. Its constructive nature allows it to generate new elements outside a given set step by step, regardless of how extensively the set is expanded.

However, as shown in this work, even an infinite iteration of such extensions does not transcend the bounds of countable cardinality. All new elements obtained diagonally belong to the constructive domain: they are enumerable, algorithmically definable, and fall within the boundaries of the set \(\mathcal{R}\)—all computable sequences.

This points to a fundamental property of the diagonal method: it expresses and proves uncountability but does not construct it. It can demonstrate the incompleteness of any algorithmic coverage but cannot overcome this incompleteness. Therefore, despite its universality as a proof technique, it cannot be used as a tool for verifying or constructing uncountable sets within constructive mathematics.

Thus, the diagonal method is a boundary between the constructive and the transcendent: it points to the existence of "beyond the enumerable" but does not allow passage into it. This also opens a direction for further analysis: comparing the expressive boundaries of constructive methods with the axiomatic strength of ZFC, where the uncountability of \(\mathbb{B}^\mathbb{N}\) is postulated but not constructively realized. This contrast may help bridge the gap between constructive logic and classical set theory.

\end{document}